\input epsf
\epsfverbosetrue

\input amssym.def
\line{}\voffset 0truein
\hoffset 0truein
\baselineskip=14pt
\def\bs{\bigskip}
\def\ms{\medskip}
\def\ss{\smallskip}
\def\c{\centerline}

\def\n{{\bf n}}
\def\ni{\noindent}

\def\<{\langle}
\def\>{\rangle}
\def\a{\alpha}
\def\b{\beta}
\def\g{\gamma}

\def \={ {\buildrel \cdot \over =}}
\def\qed{{$\vrule height4pt depth0pt width4pt$}}

\c{\bf AN INVARIANT FOR OPEN VIRTUAL STRINGS}\ss

\c{Daniel S. Silver and Susan G. Williams} \bs

\ms
\footnote{} {Both authors partially supported by NSF grant
DMS-0304971.}
\footnote{}{2000 {\it Mathematics Subject Classification.}  
Primary 57M25; Secondary 20F05, 20F34.}

{\narrower {\bf ABSTRACT:} Extended Alexander groups are used to define an invariant for open virtual   strings. Examples of non-commuting open strings and a ribbon-concordance obstruction are given. An example is given of a slice open virtual   string that is not ribbon. Definitions are extended to open $n$-strings. \bs


\ni {\bf 1. Introduction.} A {\it classical open string} is the image $\a$ of a generic immersion
${\Bbb R} \to {\Bbb R}^2$ such that $x \mapsto (x, 0)$ for $|x|$ sufficiently large. When some double points are encircled, $\a$ is said to be a {\it open virtual   string} (or simply, an {\it open string}). Two open strings are regarded as the same if they are related by a finite sequence of {\it flat generalized Reidemeister moves} (see Figure 1).  The set ${\cal M}$ of open strings is a semigroup under concatenation.\bs

\epsfxsize=3.5truein
\centerline{\epsfbox{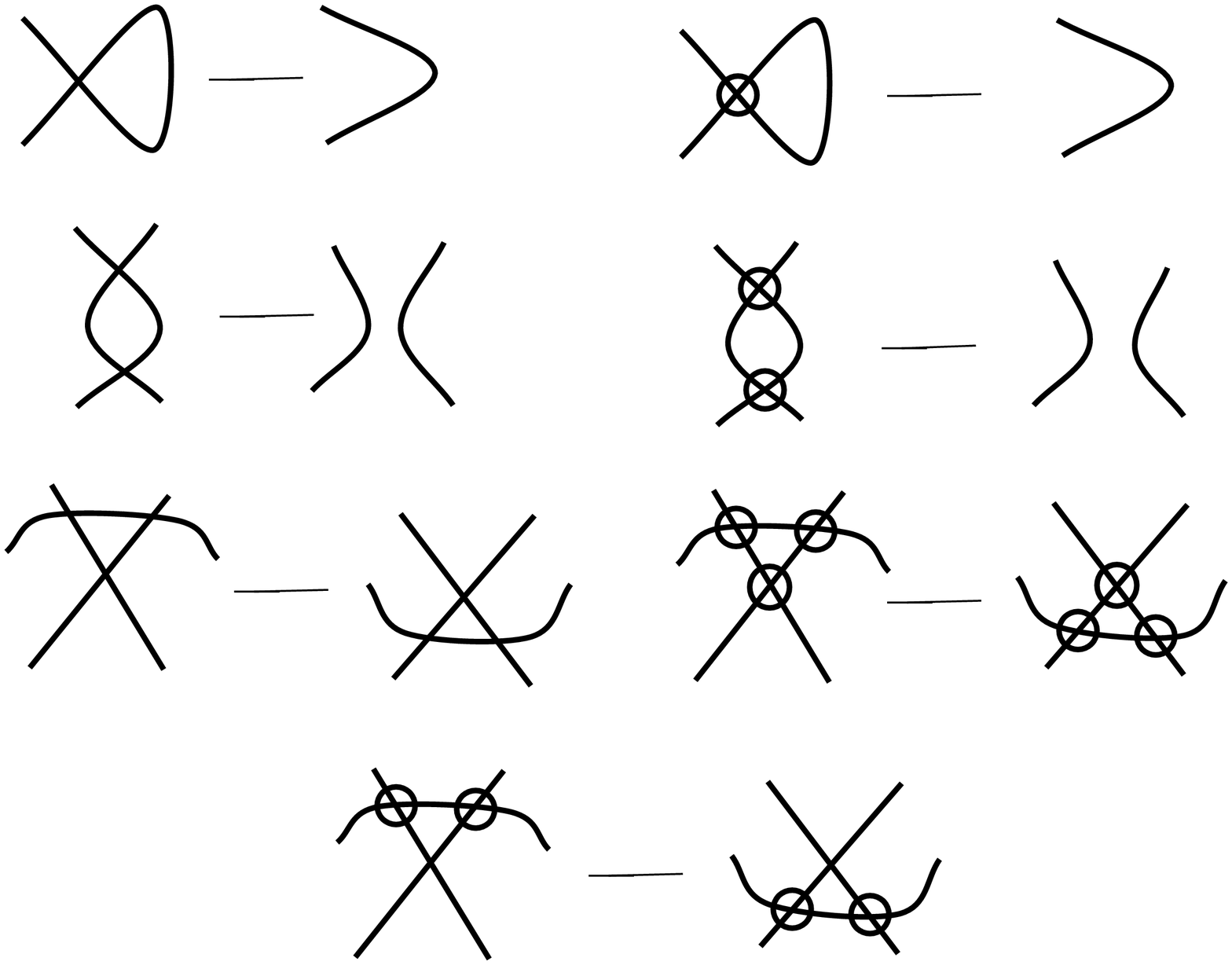}} \bs
\c{{\bf Figure 1:} Flat generalized Reidemeister moves} \bs

Gauss diagrams provide an alternative means for defining open strings. Following [T04], we describe an open string by an oriented line with $2m$ distinguished points partitioned into $m$ ordered pairs. A pair $(s,t)$ is the preimage of an uncircled double point, and the segment of {\Bbb R} near $s$ crosses the segment near $t$ transversely from left to right. We visualize each ordered pair as as an arrow above the line starting at the first coordinate of the pair and ending at the second. An equivalence relation on Gauss diagrams, called {\it homotopy}, is described in [T04], where it is shown that two open strings are the same if and only if their underlying Gauss diagrams are homotopic. \bs

Non-open strings are defined by Gauss diagrams in an analogous manner, replacing the line with an oriented circle [T04]. 
Strings and open strings represent a natural variation of the theory of virtual knots, introduced by L. Kauffman [K99].  Strings correspond to flat virtual knots while open strings can be regarded as flat virtual $1$-tangles [K99]. 
Our terminology was introduced by V. Turaev, who has borrowed ``open string" from physics and ``virtual" from Kauffman's theory.
\bs

\epsfxsize=3.75truein
\centerline{\epsfbox{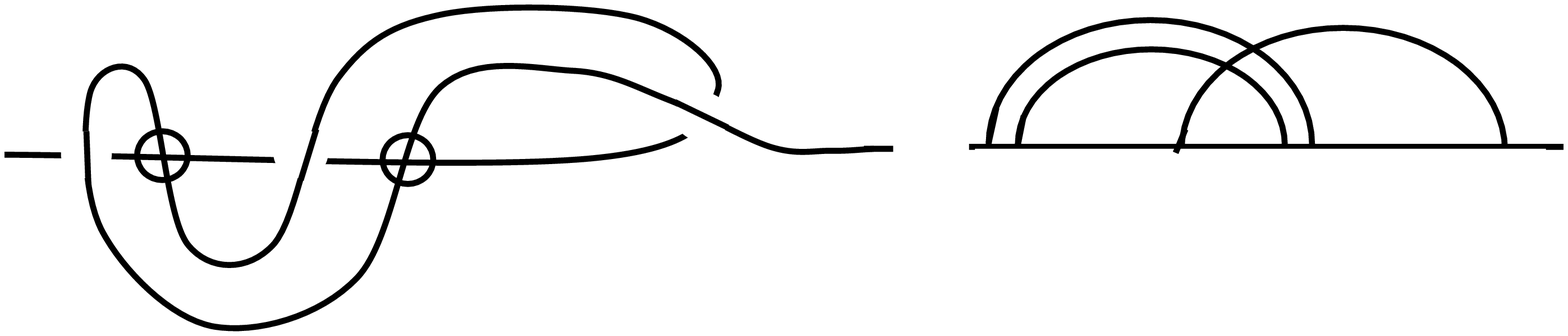}} \bs

\c{{\bf Figure 2:} A open virtual   string $\a$ and its Gauss diagram} \bs

In [SW04] we adapted Alexander group techniques of [SW01$'$] and [SW03] in order to define invariants for {\it long virtual knots}, open strings for which each uncircled double point is resolved as a classical crossing and the usual generalized Reidemeister moves for virtual knots are allowed.  Like open strings, long virtual knots form a semigroup, and the invariant of [SW04] provides the means to show that many pairs of long virtual knots do not commute. Earlier, V. Manturov [M04] had used quandle techniques to prove noncommutativity for a particular pair, and [SW04] was motivated by his work.

We modify the invariant of [SW04] to obtain a homomorphism $\Phi$ from the semigroup ${\cal M}$ of open strings to the group of endomorphisms of a countable rank free group. Using it, we show that ${\cal M}$ is non-commutative, thereby answering a question posed in [T04]. In \S3 we show that $\Phi(\a)$ has  a special form whenever $\a$ is a ribbon open knot, and we give an example of a slice open knot that isn't ribbon. In the last section, we extend the invariant $\Phi$ to open multi-strings.  

We thank Vladimir Turaev for questions, suggestions and encouragement, and also Seiichi Kamada and Naoko Kamada for valuable conversations. We are grateful to the referee for comments and suggestions that improved the paper. \bs


\ni {\bf 2. Extended Alexander groups and open strings.} The Alexander group ${\cal A}_k$ of a virtual knot was introduced in [SW01$'$]. It is a special case of R. Crowell's derived group [C84]. When $k$ is classical, ${\cal A}_k$ abelianizes to the Alexander module of the knot, a fact that provides justification for the name we have given it. The extended Alexander group $\tilde {\cal A}_k$, also defined in [SW01$'$], contains new information, and in many cases it can show that a virtual knot is not classical. 

Extended Alexander groups are defined for long virtual knots in [SW04]. 
In order to obtain an invariant of an open string $\a$, we apply the following observation of Turaev. Give $\a$ the induced orientation from ${\Bbb R}$, directing the open string from  left to right. Let $D$ be the associated descending diagram for a long virtual knot $k$, that is, while traveling along $\a$ in the preferred direction, each time an uncircled double point is encounted for the first time, resolve it as a classical crossing in such a way that we travel along the overcrossing arc.  Encircled double points are treated as virtual crossings. One easily verifies that if two open strings are related by a flat generalized Reidemeister move, then the associated long virtual knots are related by a corresponding move. Consequently, invariants of $k$ provide invariants for $\a$. 

Each overcrossing arc of the descending diagram $D$ is regarded as a union of arcs joined at the point of overcrossing. Virtual crossings are disregarded when determining arcs.  As in [SW04],  associate to $k$ a group $\tilde {\cal A}_k$. It is described by families of generators $a_\n, b_\n, c_\n, \ldots$, each indexed by vectors  $\n \in {\Bbb Z}^2$, corresponding to the arcs of $D$. To a classical positive crossing labeled as in Figure 3, associate a similarly indexed family of relations of the form 

$$a_\n b_{\n+u}= c_\n a_{\n+ u},\  c_{\n+v}=b_\n  \eqno(2.1)$$
For a negative crossing, the second relator family is replaced by $a_{\n+v}=d_\n$. Here $u, v$ are standard generators of ${\Bbb Z}^2$.

It is convenient to abbreviate generator families by $a,b,c,\ldots$, and  $a_\n b_{\n+u}= c_\n a_{\n+ u}$,  $a_{\n+v}=d_\n$ by $ab^u = ca^u$ and $a^v =d$, respectively. 
In this notation, $(j,k)\in {\Bbb Z}^2$ is written multiplicatively as $u^jv^k$.  We regard
$\tilde {\cal A}_k$ as a ${\Bbb Z}^2$-group, a group together with an action 
${\Bbb Z}^2\to {\rm Aut}(\tilde {\cal A}_k)$. \bs

\epsfxsize=2.3truein
\centerline{\epsfbox{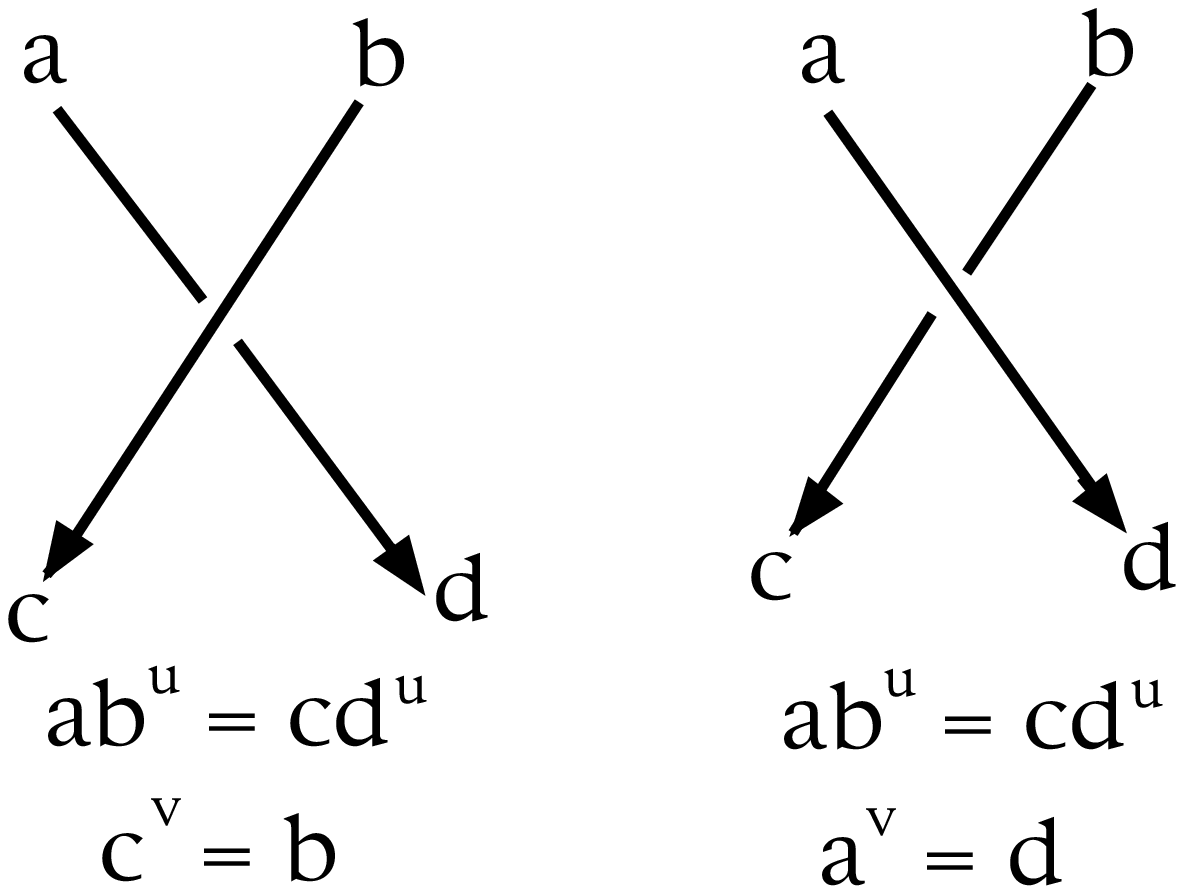}} \bs

\c{{\bf Figure 3:} Relations for extended Alexander group} \bs

Let $a_{-\infty}$ and $a_{+\infty}$ be generators corresponding to the arcs that run to $-\infty$ and $+\infty$, respectively.  We refer to these  as left and right {\it end generators}. Associate a triple $(\tilde {\cal A}_k, 
a_{-\infty}, a_{+\infty})$. Two such triples are regarded as the same if there is an isomorphism of underlying groups that commutes with the 
${\Bbb Z}^2$-actions and matches corresponding end generators. As in [SW04], the triple is an invariant of $k$. 

Since the diagram $D$ is descending, it follows easily that 
$a_{-\infty}$ freely generates $\tilde {\cal A}_k$. In particular, 
$a_{+\infty}$ is uniquely expressible as a reduced word
in the free group $F$ on generators $a_{-\infty}$.
The correspondence determines an endomorphism $\Phi(\a)$ of $F$ that respects the ${\Bbb Z}^2$-action.  Moreover, 
as in [SW04], if $\a_1$ and $\a_2$ are two open strings, then 
$\Phi(\a_1 \cdot \a_2)$ is the composition of $\Phi(\a_1)$ and $\Phi(\a_2)$. 

We summarize the above results as follows. \bs

\ni {\bf Theorem 2.1.} The assignment $\a \mapsto \Phi(\a)$ 
defines a homomorphism $\Phi: {\cal M} \to {\rm End}(F)$. \bs

\ni {\bf Remark 2.2.} $\Phi$ is similar in certain respects to the Burau representation for braids. The latter can be regarded as an automorphism of the ${\Bbb R}$-vector space spanned by the braid inputs. The image of each basis element is seen at the braid outputs. This point of view can be found, for example, in [SW01]. \bs

An open string that has no circled double points is called {\it classical}. The descending diagram $D$ will be isotopically trivial in this case. The following criterion is immediate. \bs

\ni {\bf Corollary 2.3.} If $\a$ is a classical open string, then $\Phi(\a)$ is the identity map. \bs

\ni {\bf Example 2.4} Consider the open strings $\a_1$ and $\a_2$  in Figure 4 together with the associated long virtual knots  $k_1$ and $k_2$. The knot diagrams have been labeled with generators of the extended Alexander groups. For notational convenience, we have abbreviated $a_{-\infty}$ by $a$. 

The endomorphism $\Phi(\a_1)$ maps $a$ to $\bar a^{u^{-1}v^{-2}}\bar a^{u^{-2}v^{-1}}a^{u^{-2}v^{-2}}a^{u^{-1}}a^{v^{-1}}$ while $\Phi(\a_2)$ maps $a$ to $a^v a^u a^{u^2 v^2} \bar a^{u^2 v} \bar a^{u v^2}$, where  $\ \bar{}\ $ denotes inversion. By Corollary 2.3, neither $\a_1$ nor $\a_2$ is classical.

We employ the functional notation $(x)f$ instead of $f(x)$ so that the composition 
``$f$ followed by $g$" is denoted by $fg$. It is straightforward to compute $\Phi(\a_i\a_j)=\Phi(\a_i)\Phi(\a_j)$, $i=1,2$. Simply substitute $(a)\Phi(\a_i)$ for each occurrence of $a$ in $(a)\Phi(\a_j)$. One checks that $\Phi(\a_1 \cdot \a_2)$ is not equal to $\Phi(\a_2\cdot \a_1)$. Hence the open strings $\a_1$ and $\a_2$ do not commute. \bs

\epsfxsize=4.7truein
\centerline{\epsfbox{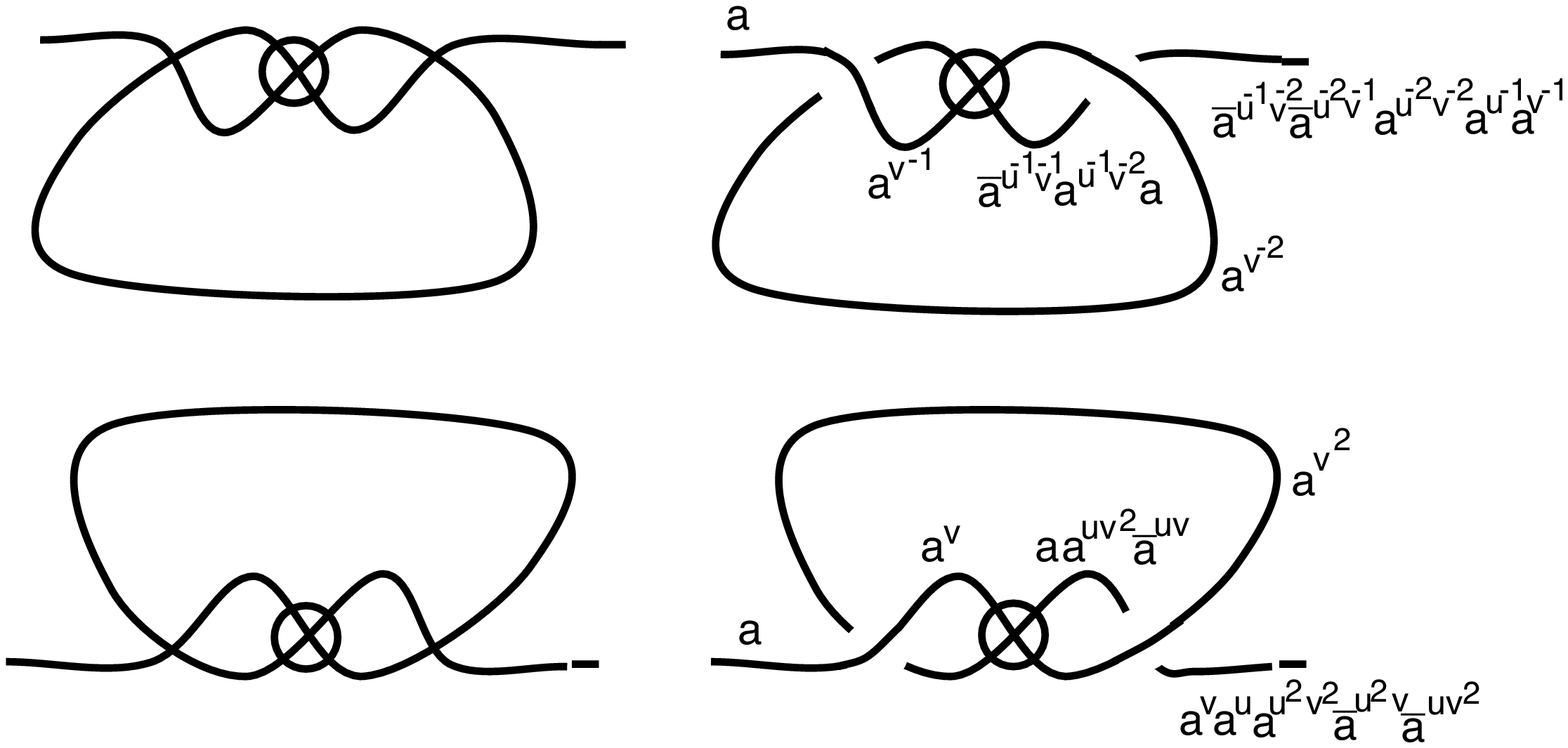}} \bs

\c{{\bf Figure 4:} Non-commuting $\a_1$ (top), $\a_2$ (bottom) and associated long virtual knots} \bs

\ni {\bf Remark 2.5.} In the above example, $a_{+\infty}$ has the form $Wa^{(uv)^{-\omega}}\bar W^{uv}$, where $W$ is a word in generators $a_{-\infty}$ and $\omega$ is equal the sum of signs of classical crossings. The fact holds generally for any open string, but we will not use it here. The proof is left to the reader.\bs

\ni {\bf Example 2.6.} The representation $\Phi$ is not faithful. Consider the open string $\b$ in Figure 5. We leave it to the reader to check that $\Phi(\b)$ is trivial. On the other hand, if we close the associated long virtual knot $k$, we obtain a virtual knot $\hat k$ that has nontrivial Jones polynomial. The latter claim can be seen easily, since $\hat k$ is  a connected sum of two simpler virtual knots, each of which has nontrivial Jones polynomial.  Hence $\a$ is a nontrivial open string. \bs

\epsfxsize=2.3truein
\centerline{\epsfbox{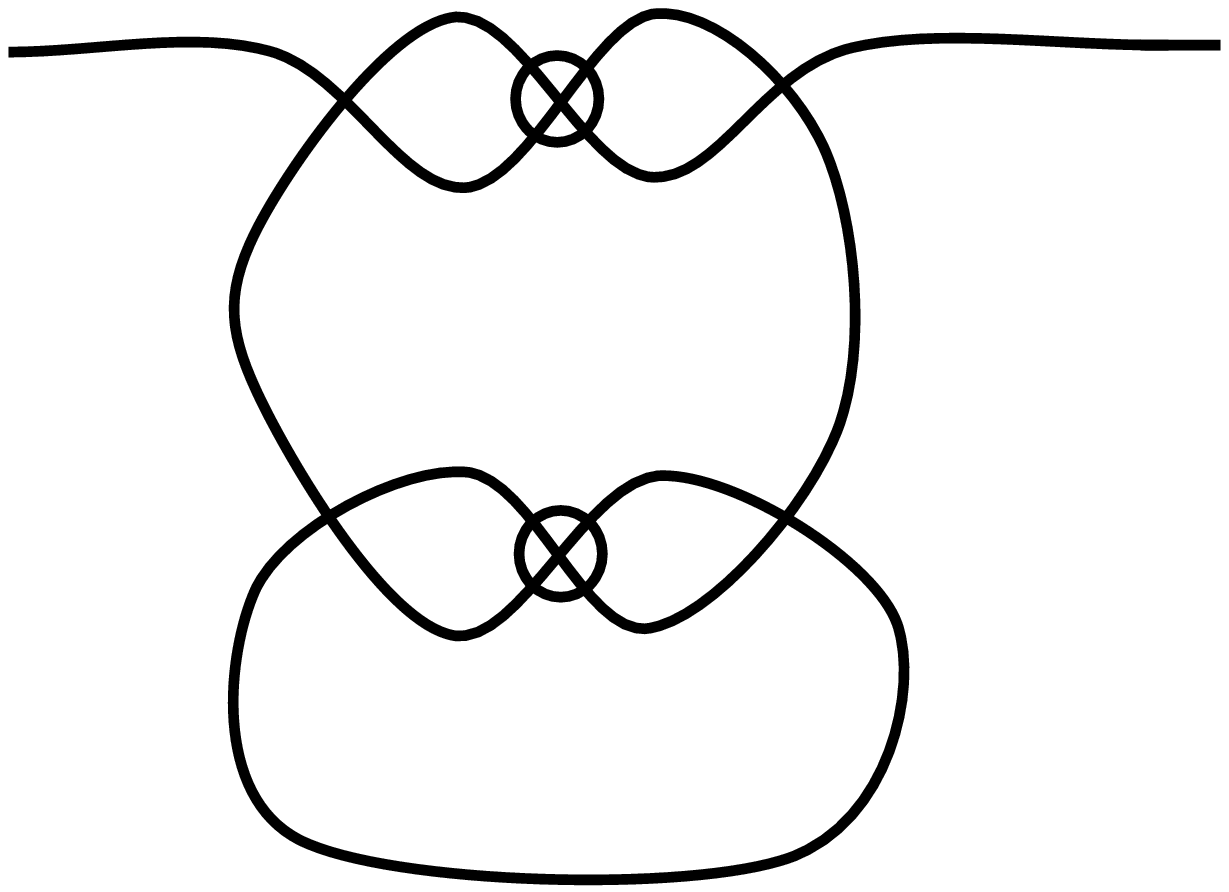}} \bs
\c{{\bf Figure 5:} Open string $\b$}\bs

\ni {\bf 3. Ribbon open strings.} As discussed in \S1, open strings can be defined either by generic immersions ${\Bbb R} \to {\Bbb R}^2$ or by Gauss diagrams. They can also be defined  in terms of generic paths in surfaces, a third perspective that motivates the notion of cobordism. Slice and ribbon open strings can be defined, and as in classical knot theory, ribbon implies slice [T04]. \bs

\ni {\bf Definition 3.1.} An open string is {\it ribbon} if it has a Gauss diagram for which reflection about the origin takes each arrow to an arrow with opposite direction. \bs

The polynomial invariant in \S3.2 of [T04] provides an elementary obstruction for an open string to be slice (and hence an obstruction for it to be ribbon). It vanishes for the non-ribbon open string $\a$ in the following example.\bs

\ni {\bf Example 3.2.} Consider the ribbon open strings $\a_2$ and $\a_3$ and their Gauss diagrams in Figure 6. (Here $\a_2$ is the same as $\a_2$ in Example 2.4.) Let $\a$ be the product $\a_2 \cdot \a_3$. A Gauss diagram for $\a$ is obtained by concatenating the Gauss diagrams for $\a_2$ and $\a_3$. Reflecting the Gauss diagram for $\a$ and changing the direction of each arrow produces a Gauss diagram for $\a_3 \cdot \a_2$. Hence if $\a$ is ribbon, then 
$\a_2 \cdot \a_3$ is equal to $\a_3 \cdot \a_2$.

Computing $\Phi(\a_3)$ as in Example 2.4, and using similar notation, we find that
$(a)\Phi(\a_3) = a a^{uv^2} a^{u^2v} a^{u^3v^3}\bar a^{u^3v^2}\bar a^{u^2v^3}\bar a^{uv}.$ The reader can check that 
$\Phi(\a_2)$ and $\Phi(\a_3)$ are non-commuting endomorphisms. Hence $\a_2$ and $\a_3$ do not commute, and so $\a$ is not ribbon. 

Since ribbon implies slice for open strings, and since a product of slice open strings is slice, $\a$ is an example of a slice open string that is not ribbon. \bs

\epsfxsize=5truein
\centerline{\epsfbox{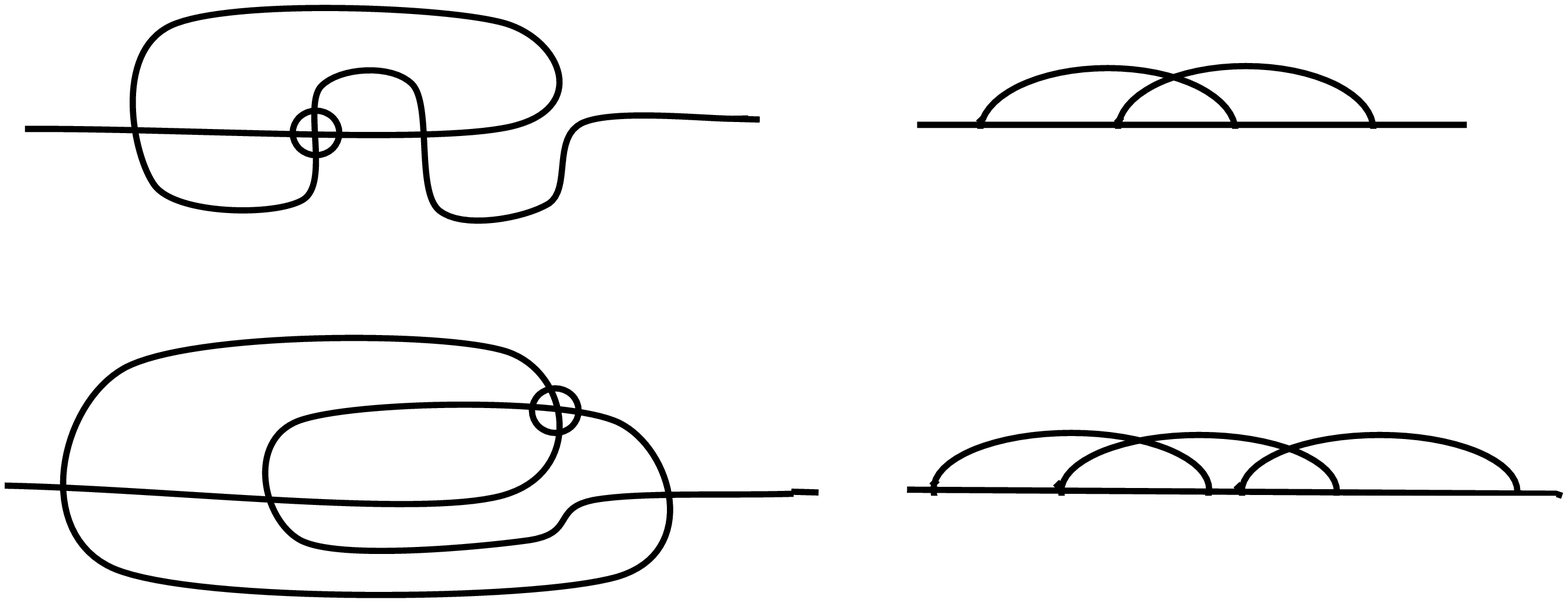}} \bs
\c{{\bf Figure 6:} Open strings $\a_2$ (top) and $\a_3$ (bottom)}\bs

An obstruction for a general open string $\a$ to be ribbon can be described. First compute $\Phi(\a)$. Then reflect the Gauss diagram for $\a$ about the $y$-axis and change the direction of each arrow, thereby obtaining a second open string $\hat\a$. If its image under $\Phi$ does not agree with $\Phi(\a)$, then $\a$ cannot be ribbon.

Testing to see if a given open string is ribbon requires two, sometimes tedious calculations. However, by considering the abelianized endomorphism $\Phi$, we obtain an obstruction that requires only one, relatively easy calculation, and yet remains very effective. In order to describe it, we let $\phi(\a)(u,v)\in {\Bbb Z}[u^{\pm 1}, v^{\pm 1}]$ be the Laurent polynomial given as the sum of the exponents of $(a)\Phi(\a)$. 
For example, the ribbon open string $\a_2$ in Example 2.4 has polynomial $v+u+u^2v^2 -u^2v-uv^2$. \bs

\ni {\bf Proposition 3.3.} If $\a$ is an open ribbon string, then $\phi(\a)(u,v) = \phi(\a)(v,u)$. \bs

\ni {\bf Proof.} We first describe a simple way to calculate $\phi(\a)$ from the Gauss diagram. Let $D$ be the corresponding descending diagram of a long virtual knot $k$. Write the generators $\tilde {\cal A}_k$ corresponding to arcs of $D$ as words in $\langle \a_{-\infty}\rangle$, as in Figure 4, and associate to each arc $x$ the Laurent polynomial $p_x(u,v)$ given by the exponent sum, so that $\phi(\a)= p_{a_{+\infty}}$. Solving the crossing relations in Figure 3 for $c$ and $d$ and abelianizing shows that $p_c=v^{-1}p_b$, $p_d=u^{-1} p_a + (1-u^{-1}v^{-1})p_b$ for the positive crossing and $p_d= v p_a$, $p_c=up_b+(1-uv)p_a$ for the negative one. 

Positve and negative crossings in $D$ correspond to arrows $(s,t)$ in the Gauss diagram of $\a$ with $s<t$ and $s>t$, respectively. We indicate the crossing relations on the Gauss diagram on the left in Figure 7. Place weights on the arrows and adjoining line segments as shown on the right in Figure 7. (A segment may receive zero, one or two weights as in Figure 8.) Consider all paths that traverse the Gauss diagram from left to right following some combination of arrows and line segments. For each path, we compute the product of the weights along the path. Then $\phi(\a)$ is the sum of the these products over all paths. We can see this by observing inductively that each $p_x$  is  determined in this way by the paths leading to $x$.

If we reflect the weighted Gauss diagram about the $y$-axis, reverse all arrows and then replace $u$ and $v$ by $v$ and $u$, respectively, we obtain the weighted Gauss diagram for $\hat \a$. Thus $\phi(\hat \a)(u,v) = \phi(\a)(v,u)$, and this is equal to $\phi(\a)(u,v)$ if $\a$ is ribbon. \qed\bs

\epsfxsize=5truein
\centerline{\epsfbox{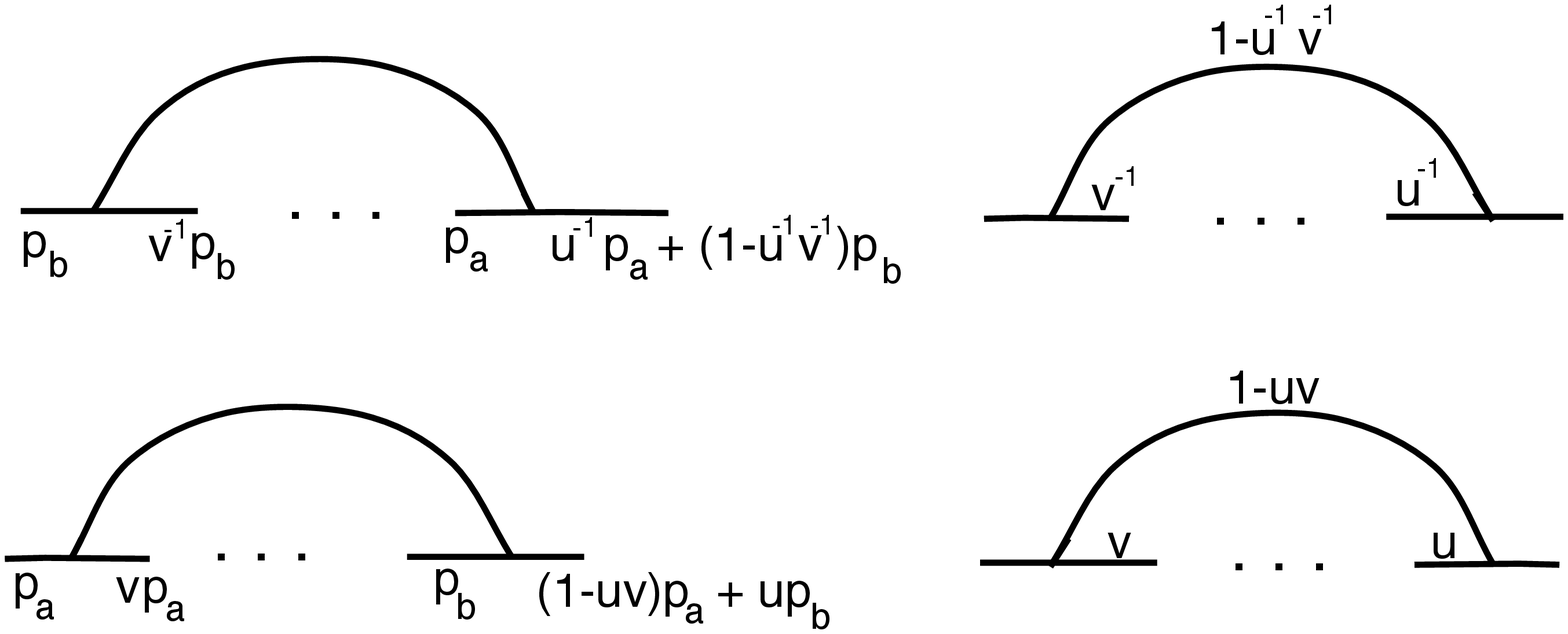}} \bs
\c{{\bf Figure 7:} Weights for Gauss diagrams}\bs

\ni {\bf Example 3.4.} For the open string $\a$ with weighted Gauss diagram in Figure 8 there are five paths traversing the diagram, and $\phi(\a)(u,v)$ is given by
$$vvuv^{-1}uu^{-1}+ (1-uv)v^{-1}uu^{-1}+v(1-uv)u^{-1}+
vvu(1-u^{-1}v^{-1})+(1-uv)(1-u^{-1}v^{-1})$$
$$=uv^2-v^2-u-v+u^{-1}v -v^{-1}.$$
Since $\phi(\a)(u,v)\ne \phi(\a)(v,u)$, the open string $\a$ is not ribbon. \bs

\epsfxsize=3.5truein
\centerline{\epsfbox{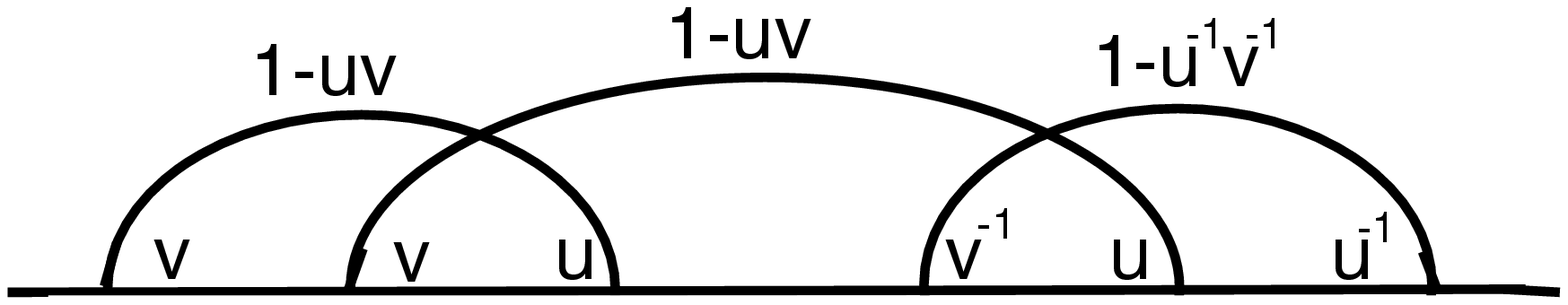}} 

\c{{\bf Figure 8:} Weighted segments of a Gauss diagram}\bs

Changing the directions of all arrows of a Gauss diagram for an open string $\a$ is equivalent to reflecting $\a$ about a line in the plane. (This is easily seen by choosing the line to be the $x$-axis.) We denote the reflected open string by $\a*$. The descending diagrams $D, D^*$ for the associated long virtual knots are also reflections of each other, and the effect of reflection on the invariant $\Phi$ is easily described. If $W$ is a word in the generators $a_\n = a^{u^jv^k}$, let $W^*$ be the word obtained by writing $W$ in reverse order and then replacing $u,v$ by $u^{-1}, v^{-1}$, respectively. \bs

\ni {\bf Lemma 3.5.} $(a)\Phi(\a*) = [(a)\Phi(\a)]^*$. \bs

\ni {\bf Proof.} The arc labeling of $D^*$ is obtained from that of $D$ by the involution $x \mapsto x^*$. To see this, it suffices to check that if the arc labelings in Figure 9a satisfy the crossing relation (2.1), then the labelings in the reflected diagram of Figure 9b satisfy the corresponding crossing relation for negative crossings. The calculation is routine. \qed\bs

\epsfxsize=2.5truein
\centerline{\epsfbox{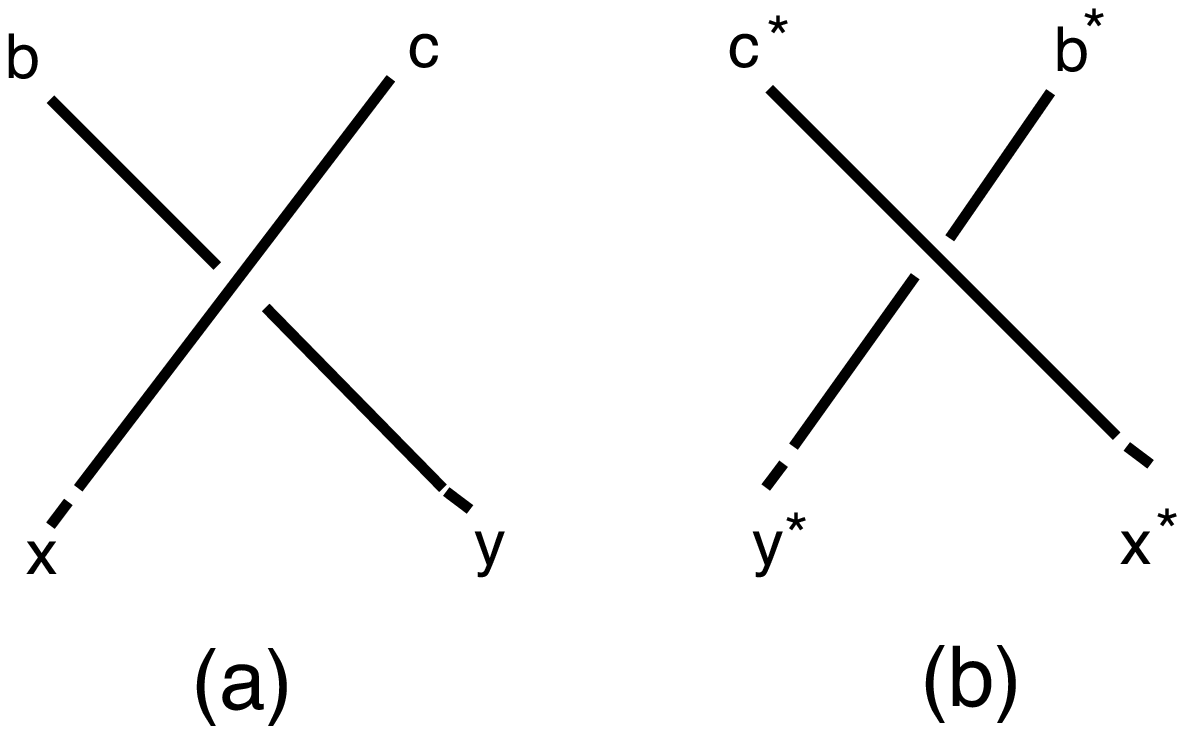}} \bs
\c{{\bf Figure 9:} Labeled crossing diagram and its reflection}\bs

\ni {\bf Example 3.6.} Consider the open string $\a$ described by the Gauss diagram in Figure 2. Let $\b$ be the product $\a\cdot \a^*$. As in Example 3.2, the slice obstruction in [T04] vanishes. We apply Proposition 3.3 to show that $\b$ is not ribbon.

We have $\phi(\b)(u,v) = \phi(\a)(u,v) \cdot \phi(\a^*)(u,v) = \phi(\a)(u,v) \cdot \phi(\a)(u^{-1}, v^{-1})$. One easily computes that $\phi(\a)(u,v)= -uv^3 - u^3v^2 +u^3v^3+u+v^2$, and from this that $\phi(\b)(u,v)$ is not equal to $\phi(\b)(v,u)$. Hence $\g$ is not ribbon. \bs

\ni {\bf Remark 3.7.} The abelian invariant $\phi$ cannot detect non-commutativity of open strings, and so it is not effective for the construction in Example 3.2.\bs


\ni {\bf 4. Open virtual   multi-strings.} Open virtual   strings generalize in several ways. For example, we can consider several open components.  

An {\it open virtual   $n$-string} (or simply {\it open $n$-string}) is the image $\a$ of a generic immersion ${\Bbb R}\times \{1, \ldots, n\} \to {\Bbb R}^2$ such that $(x, j) \mapsto (x, j)$ for sufficiently large negative $x$, and $(x, j) \mapsto (x, \pi_j)$ for sufficiently large positive $x$. Here $(\pi_1, \dots, \pi_n)$ is a permutation of $(1,\ldots, n)$. Some double points may be encircled, and as in the special case of open strings, when $n=1$, two open $n$-strings are regarded as the same if they are related by a finite sequence of flat generalized Reidemeister moves (Figure 1).

In classical link theory, one often puts labels on components and keeps track of them. 
Similarly, we assign integers $1,2,\ldots, n$ (called {\it colors}) to the components of an open $n$-string $\a$, and we refer to $\a$ as a {\it colored open $n$-string}.  The component of $\a$ labeled $i$ is the {\it $i$th component}. {\sl We multiply two colored open $n$-strings  only if the colors of each pair of joined components agree.} 

To each colored open $n$-string $\a$ we assign a ${\Bbb Z}^{n+1}$-group in a manner similar to that of \S2. Begin by constructing a descending diagram $D$ for a long virtual link $\ell$: While traveling along the first component of $\a$ in the preferred direction,  each time an uncircled double point of $\a$ is encounted for the first time, resolve it as a classical crossing in such a way that we travel along the overcrossing arc. Repeat with the remaining components in turn, similarly resolving the uncircled double points that remain. Double points are treated as virtual crossings. As before, one verifies that if two colored open $n$-strings are related by a flat generalized Reidemeister move, then the associated long virtual links are related by a corresponding move. Again, invariants of $\ell$ provide invariants of $\a$.

Now define an extended Alexander group $\tilde {\cal A}_k$ by associating ${\Bbb Z}^{n+1}$-group generators to arcs of $D$. At each classical crossing, associate relators as in Figure 7. \bs

\epsfxsize=3truein
\centerline{\epsfbox{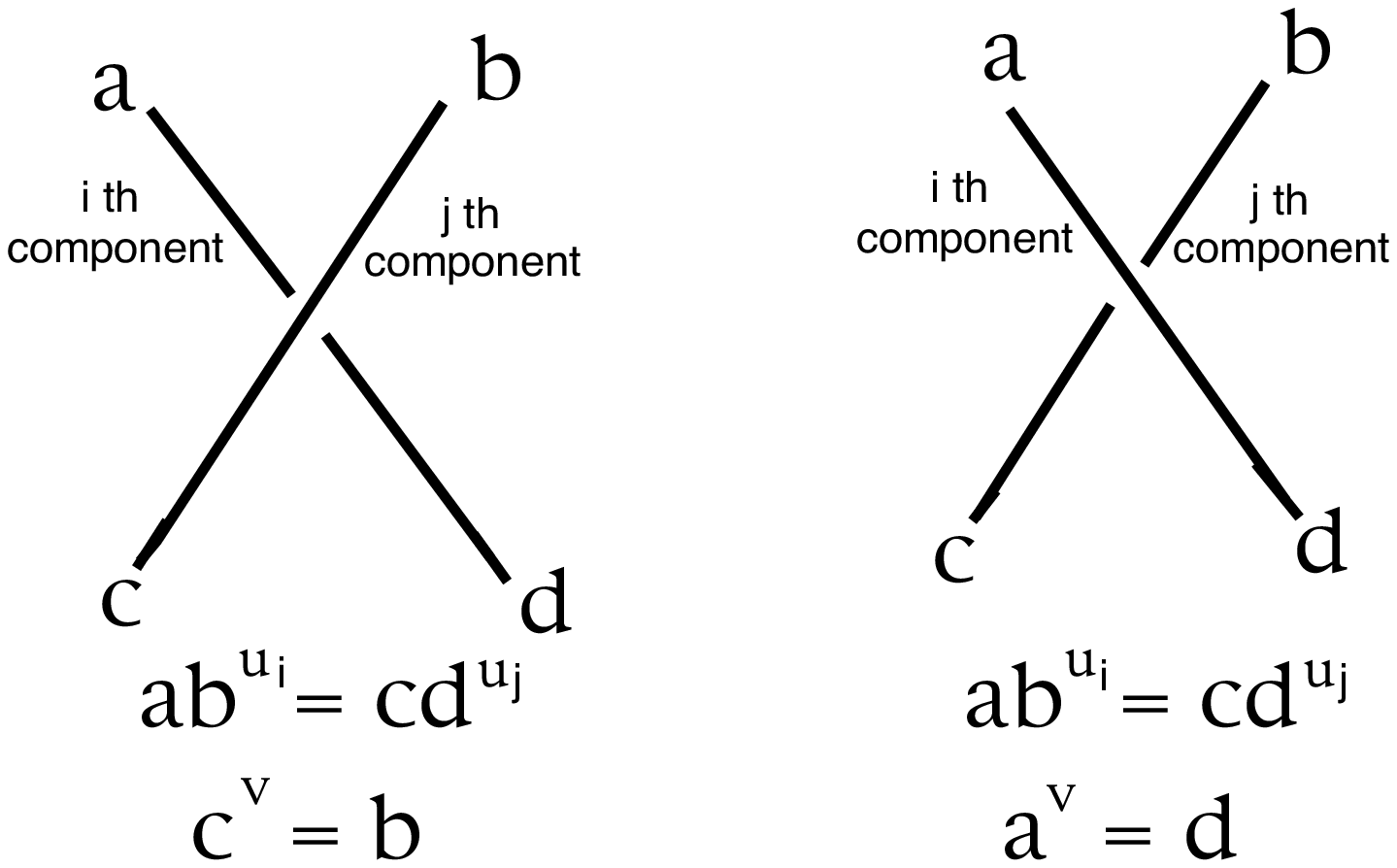}} \bs
\c{{\bf Figure 10:} Extended Alexander group relations for an open $n$-string} \bs

Let $a_{-\infty}$ and $a_{+\infty}$ be generators corresponding to the arcs on the first component of $D$ that  are leftmost and rightmost, respectively. Let 
$b_{-\infty}$ and $b_{+\infty}$ be similar generators for the second component, and so forth. We will refer to these  as left and right {\it end generators}. Associate a tuple $(\tilde {\cal A}_k, 
a_{-\infty}, b_{-\infty},\ldots, a_{+\infty}, b_{+\infty}, \ldots )$. Regard two such tuples as the same if there is an isomorphism of underlying groups that commutes with the 
${\Bbb Z}^{n+1}$-actions and matches corresponding end generators. As in [SW04], the tuple is an invariant of $\ell$. (This is proven only in the case that $\ell$ has a single component. However, the general argument is essentially the same.) 

By design, the diagram $D$ is descending. It follows as before that 
$a_{-\infty}, b_{-\infty},\ldots$ generate $\tilde {\cal A}_k$. In particular, 
$a_{+\infty}, b_{+\infty}, \ldots$ are uniquely expressible as words
in the free group $F$ on $a_{-\infty}, b_{-\infty},\ldots$.
Again, the correspondence determines an endomorphism $\Phi(\a)$ of $F$.  Moreover,  if $\a_1$ and $\a_2$ are two colored open $n$-strings such that the product $\a_1\cdot \a_2$ is defined, then 
$\Phi(\a_1 \cdot \a_2)$ is the composition $\Phi(\a_1) \circ \Phi(\a_2)$.

\bs
\ni {\bf References.} \bs

\ni [C84] R. Crowell, The derived group of a permutation representation, {\sl Advances in Math.\ \bf53} (1984), 88--124. \ss


\ni [K99] L.H. Kauffman, Virtual knot theory, {\sl European J. Comb.\ \bf 20} (1999), 663--690. \ss

\ni [M04] V. Manturov, Long virtual knots and their invariants II: the commutation problem, preprint. \ss


\ni [SW01] D.S. Silver and S.G. Williams, A generalized Burau representation for open string links, {\sl Pacific J. Math.\ \bf197} (2001), 241--255.\ss

\ni [SW01$'$] D.S. Silver and S.G. Williams, Alexander groups and virtual links, {\sl J. Knot Theory and its Ramifications\ \bf 10} (2001), 151--160. \ss

\ni [SW03] D.S. Silver and S.G. Williams, Polynomial invariants of virtual links, {\sl J. Knot Theory and its Ramifications\ \bf 12} (2003), 987--1000. \ss

\ni [SW04] D.S. Silver and S.G. Williams, Alexander groups of long virtual knots, {\sl J. Knot Theory and its Ramifications}, to appear; arXiv:math.GT/0405460\ss

\ni [T04] V. Turaev, Virtual  strings and their cobordisms, preprint, 2004,  arXiv: math.GT/\hfil\break0311185    \bs

\ni address: Dept. of Mathematics and Statistics, Univ. of South Alabama, Mobile AL 36688 USA \ss

\ni e-mail: silver@jaguar1.usouthal.edu; swilliam@jaguar1.usouthal.edu

\end